\documentclass[11pt]{amsart}

\usepackage{fullpage}
\usepackage[utf8]{inputenc}
\usepackage[T1]{fontenc}
\usepackage{graphicx}
\usepackage{amsmath,amsfonts,amssymb,mathtools}
\usepackage{stmaryrd}
\usepackage{graphicx}
\usepackage{mathrsfs,dsfont}
\usepackage{amsmath,amsthm,amsthm}

\usepackage{hyperref}

\usepackage{enumerate}

\usepackage{color}
\usepackage{amsmath,amsfonts,amssymb,mathtools}
\usepackage{stmaryrd}
\usepackage{graphicx}
\usepackage{mathrsfs,dsfont}
\usepackage{amsmath,amsthm,amsthm}

\usepackage{hyperref}

\newcommand{\nrep}{n_{\rm rep}}

\newtheorem{theo}{Theorem}

\newtheorem{rem}{Remark}

\begin{document}


\title{On a new class of score functions to estimate tail
  probabilities of some stochastic processes
with Adaptive Multilevel Splitting}



\author{Charles-Edouard Br\'ehier}
\address{Univ Lyon, CNRS, Université Claude Bernard Lyon 1, UMR5208, Institut Camille Jordan, F-69622 Villeurbanne, France}
\email{brehier@math.univ-lyon1.fr}

\author{Tony Leli\`evre}
\address{Universite\'e Paris-Est, CERMICS (ENPC), INRIA, 6-8-10 avenue Blaise Pascal, 77455 Marne-la-Valle\'ee,
France}
\email{lelievre@cermics.enpc.fr}

\date{\today}

\begin{abstract}
We investigate the application of the Adaptive Multilevel Splitting algorithm for the estimation of tail probabilities of solutions of Stochastic Differential Equations evaluated at a given time, and of associated temporal averages.

We introduce a new, very general and effective family of score functions which is designed for these problems. We illustrate its behavior on a series of numerical experiments. In particular, we demonstrate how it can be used to estimate large deviation rate functionals for the longtime limit of temporal averages.
\end{abstract}


\maketitle 

\section{Introduction}

Fast and accurate estimation of rare event probabilities, and the
effective simulation of these events, is a challenging computational
issue, which appears in many fields of science and engineering. Since
rare events are often the ones that matter in complex systems, designing efficient and easily implementable algorithms is a crucial question which has been the subject of many studies in the recent years.

Since the pionnering works on Monte-Carlo methods, several classes of
algorithms have been developed, see for instance the
monographs~\cite{AsmussenGlynn2007,Bucklew2004,RubinoTuffin2009}. The
most popular strategies are importance sampling and splitting. On the
one hand, importance sampling consists in changing the probability
distribution, such that under the new probability distribution the
events of interest are not rare anymore. Appropriate reweighting then
yields consistent estimators. This strategy has for instance been
applied recently to simulate rare events in climate
models~\cite{RagoneWoutersBouchet:18}. On the other hand, splitting
techniques, consist in writing the rare event probability as a product
of conditional probabilities which are simpler to estimate, and in
using interacting particle systems in order to estimate these
conditional probabilities.

In this manuscript, a class of splitting algorithms is considered. Splitting techniques have been introduced in the 1950s~\cite{KahnHarris1951}, and have been studied extensively in the last two decades~\cite{CerouDelmoralFuronGuyader:12,DeanDupuis2009,GarvelsKroeseVan-Ommeren2002,GlassermanHeidelbergerShahabuddinZajic1999}. Many variants have appeared in the literature: Generalized multilevel splitting~\cite{BotevKroese2012,BotevLecuyerRubinoSimardTuffin:13}, RESTART~\cite{Villen-Altamirano1991,Villen-Altamirano1994}, Subset simulation~\cite{AuBeck2001}, Nested sampling~\cite{Skilling2006,Skilling2007}, Reversible shaking transformations with interacting particle systems \cite{AgarwalDeMarcoGObetLiu:18,GobetLiu:15}, genealogical particle analysis~\cite{Del-MoralGarnier2005,WoutersBouchet:16}, etc...

The Adaptive Multilevel Splitting (AMS)
algorithm~\cite{CerouGuyader:07} is designed to estimate rare event
probabilities of the type $\mathbb{P}(\tau_B<\tau_A)$, where $\tau_A$
and $\tau_B$ are stopping times associated with a Markov process $X$,
typically  the entrance times of $X$ in regions $A$ and $B$ of the state space. In many applications, $A$ and $B$ are metastable states for the process. The algorithm is based on selection and mutation mechanisms, which leads to the evolution of a system of interacting replicas. The selection is performed using a score function, which is often referred to as a reaction coordinate when dealing with metastable systems.

The objective of this article is to design and test new score
functions, using the AMS strategy, to estimate probabilities of the
type $\mathbb{P}(\Phi(X_T)>a)$ or
$\mathbb{P}\left(\frac{1}{T}\int_{0}^{T}\phi(X_t)dt>a\right)$, where
$a$ is a threshold, $\Phi$ and $\phi$ are real-valued functions, and
$\bigl(X_t\bigr)_{0\le t\le T}$ is a Markov process. In fact, as will
be explained below, the probabilities of interest can be rewritten as
$\mathbb{P}(\tau_B<\tau_A)$, associated with an auxiliary Markov
process. Our main contribution is the identification of appropriate
score functions related to this interpretation, and which return a
non-zero value for the estimator of the probability of interest. By
using then an AMS algorithm which fits in the Generalized Adaptive
Multilevel Splitting
framework~\cite{BrehierGazeauGoudenegeLelievreRousset:16}, one can construct unbiased estimators of the probability (and possibly of other quantities of interest).

The efficiency of the approach is investigated with numerical
experiments, using several test cases taken from the literature on
rare events. First, validation is performed on one-dimensional
Gaussian models (Brownian Motion~\cite{Vanden-EinjdenWeare2012},
Ornstein-Uhlenbeck process~\cite{WoutersBouchet:16}). More complex
test cases then illustrate the efficiency of the approach and of the
new score functions introduced in this article: drifted Brownian
Motion, three-dimensional Lorenz
model~\cite{BeckZuev:17,Lorenz:63}. Estimations of probabilities
depending on temporal averages $\frac1T\int_0^T\phi(X_s)ds$ are also
considered for two models~\cite{FerreTouchette:18}: the
one-dimensional Ornstein-Uhlenbeck process and a driven periodic
diffusion. In these examples, values of large deviations rate
functionals for the longtime limit $T \to \infty$ are estimated.

In the last decade, many works have been devoted to the analysis and
applications of AMS algorithms. A series of work has been devoted to
the analysis of the so-called ideal
case~\cite{BrehierLelievreRousset:15,BrehierGoudenegeTudela:16,Brehier:16,GuyaderHentgartnerMatzner-Lober:11},
namely when the AMS algorithm is applied with the optimal score
function (namely the so-called committor function). In practice, this
optimal score function is unkown. Beyond the ideal case,
consistency~\cite{BrehierGazeauGoudenegeLelievreRousset:16}
(unbiasedness of the small probability estimator) and
efficiency~\cite{CerouDelyonGuyaderRousset:18} (variance of the small
probability estimator) have been studied. Moreover, the adaption of
the original algorithm to the discrete-in-time setting has been
studied in details
in~\cite{BrehierGazeauGoudenegeLelievreRousset:16}. It can be used
to compute transition times between metastable states~\cite{CerouGuyaderLelievrePommier:11}, return times~\cite{LestangRagoneBrehierHerbertBouchet:18}, or other
observables associated with the rare event of
interest~\cite{Louvin:17}.  The AMS algorithm has been successfully
applied in many contexts: the Allen-Cahn stochastic partial
differential
equation~\cite{BrehierGazeauGoudenegeRousset:15,RollandBouchetSimonnet:16},
the simulation of Bose-Einstein condensates~\cite{Poncet:17},
molecular dynamics and computational
chemistry~\cite{AristoffLelievreMayne:15,CerouGuyaderLelievrePommier:11,LopesLelievre:17,TeoMayneSchultenLelievre:16},
nuclear
physics~\cite{LouvinDumonteilLelievreRoussetDiop:17,LouvinDumonteilLelievre:17,Louvin:17}
and turbulence~\cite{BouchetRollandSimonnet:18,Rolland:18}, for example.

This article is organized as follows. Section~\ref{sec:setting} presents the precise mathematical setting, in particular the rare event probability of interest is defined by~\eqref{eq:proba3}. A general formulation of the AMS algorithm designed to estimate this quantity is provided in Section~\ref{sec:generalAMS}, in particular see Section~\ref{sec:algo} for the full algorithmic description. Examples of appropriate score functions are discussed in Section~\ref{sec:choice}. To overcome the limitations of a vanilla strategy, Section~\ref{sec:choice_vanilla}, our main contribution is the construction of the score functions presented in Section~\ref{sec:choice_new}. Finally, numerical experiments are reported in Section~\ref{sec:numerical}.

\section{Setting}\label{sec:setting}

We consider stochastic processes, with values in $\mathbb{R}^d$, in
dimension $d\in\mathbb{N}$, which are solutions of Stochastic
Differential Equations (SDEs) of the type: for $0 \le t_0 \le t \le T$
and $x_0 \in \mathbb R^d$,
\begin{equation}\label{eq:SDE}
dX_t^{t_0,x_0}=f(t,X_t^{t_0,x_0})dt+\sigma(t,X_t^{t_0,x_0})dW(t)
\end{equation}
where $X_t^{t_0,x_0} \in \mathbb{R}^d$, with initial condition,
\begin{equation}\label{eq:SDE-IC}
X_{t_0}^{t_0,x_0}=x_0.
\end{equation}
The noise $\bigl(W(t)\bigr)_{t\ge 0}$ is given by a standard Wiener
process with values in $\mathbb{R}^D$, for some $D\in\mathbb{N}$. The
coefficients $f:[0,T]\times \mathbb{R}^d\to \mathbb{R}^d$ and
$\sigma:[0,T]\times \mathbb{R}^d\to \mathbb{R}^{d \times D}$ are assumed to be sufficiently smooth to ensure global well-posedness of the SDE.

In this work, two types of rare events associated with
$(X_{t}^{t_0,x_0})_{0 \le t \le T}$ are considered. Let $a\in\mathbb{R}$ denote a threshold, and let $\Phi,\phi:\mathbb{R}^d\to\mathbb{R}$ be two measurable functions. First, we are interested in tail probabilities for the random variable $\Phi(X_T^{t_0,x_0})$, namely in
\begin{equation}\label{eq:proba1}
\mathbb{P}\bigl(\Phi(X_T^{t_0,x_0})>a\bigr).
\end{equation}
Second, we are interested in tail probabilities for temporal averages, defined as
\begin{equation}\label{eq:proba2}
\mathbb{P}\left(\frac{1}{T-t_0}\int_{t_0}^{T}\phi(X_t^{t_0,x_0})dt>a\right).
\end{equation}
We will investigate numerically the performance of AMS estimators for
both~\eqref{eq:proba1} and~\eqref{eq:proba2} on various examples. In
particular, we will consider the regime $T \to \infty$
for~\eqref{eq:proba2} in order to estimate large deviation rate functionals.

Notice that the case of temporal averages~\eqref{eq:proba2} can be
rewritten in the form of~\eqref{eq:proba1}. Indeed, the probability~\eqref{eq:proba2} may be written as~\eqref{eq:proba1} for the auxiliary process defined by $\tilde{X}_t^{t_0,x_0}=\bigl(X_t^{t_0,x_0},Y_{t}^{t_0,x_0}\bigr)$, where $Y_t^{t_0,x_0}=\phi(x_0)$, and for $t>t_0$
\[
Y_t^{t_0,x_0}=\frac{1}{t-t_0}\int_{t_0}^{t}\phi(X_s^{t_0,x_0})ds,
\] 
and with $\tilde{\Phi}(x,y)=y$. The process $\bigl(Y_t^{t_0,x_0}\bigr)_{t_0\le t\le T}$ is solution of the following ODE,
\[
dY_t^{t_0,x_0}=\frac{1}{t-t_0}\bigl(\phi(X_t^{t_0,x_0})-Y_t^{t_0,x_0}\bigr),~t>t_0\quad,~Y_{t_0}^{t_0,x_0}=\phi(x_0),
\]
coupled with the SDE for the diffusion process
$\bigl(X_t^{t_0,x_0}\bigr)_{t_0\le t\le T}$. This trick will be used
for our numerical experiments below. Therefore, in the following, we present the AMS algorithm and discuss its theoretical properties only for~\eqref{eq:proba1}.

For future purposes, observe that the target probability~\eqref{eq:proba1} may be written as $u_a(t_0,x_0)$, where
\begin{equation}\label{eq:Kolmogorov}
u_a(t,x)=\mathbb{P}\bigl(\Phi(X_T^{t,x})>a\bigr)
\end{equation}
is the solution (under appropriate regularity assumptions) of the backward Kolmogorov equation
\begin{equation}\label{eq:KB}
\begin{cases}
\frac{\partial u_a(t,x)}{\partial t}+\mathcal{L}_t u_a(t,x)=0 \text{
  for } t \in [0,T] \text{ and } x \in \mathbb R^d,\\
u_a(T,x)=\mathds{1}_{\Phi(x)>a} \text{ for } x \in \mathbb R^d,
\end{cases}
\end{equation}
where the infinitesimal generator $\mathcal{L}_t$ is defined by: for
all test functions $\varphi$, $\mathcal{L}_t\varphi(x)=f(t,x)\cdot\nabla
\varphi(x)+\frac{1}{2}\sigma(t,x)\sigma(t,x)^\star :
\nabla^2\varphi(x)$. Approximating the solutions of PDEs of this type
using deterministic methods is in general possible only when the
dimension $d$ is small. Instead, Monte Carlo methods may be
used. However, naive Monte Carlo algorithms are not efficient in the rare event regime, {\it e.g.} when $a\to \infty$ or when the diffusion coefficient is of the type $\sigma_\epsilon=\sqrt{\epsilon}\sigma$ and $\epsilon\to 0$.

\bigskip

In practice, discrete-time approximations are implemented. Let $\Delta
t>0$ denote the time-step size of the integrator (for instance the
standard Euler-Maruyama method), with $T=N\Delta t$ and $t_0=n_0\Delta
t$, where $n_0\in\mathbb{N}_0,N\in\mathbb{N}$, $n_0\le N$. With a
slight abuse of notation, let us denote the discrete-time process
obtained after discretization of~\eqref{eq:SDE} by
$\bigl(X_n^{n_0,x_0}\bigr)_{n_0\le n\le N}$. The time-discrete
counterpart of \eqref{eq:proba1} is then
\begin{equation}\label{eq:proba3}
\mathbb{P}\bigl(\Phi(X_N^{n_0,x_0})>a\bigr).
\end{equation}
The algorithms presented below are used to estimate probabilities of the type~\eqref{eq:proba3}.

\begin{rem}
It is assumed that the initial condition is deterministic:
$X^{n_0,x_0}_{n_0}=x_0$. The adaptation of the algorithms presented
below to the case of a random initial condition is straightforward, by
simply using the Markov property:
$\displaystyle{\mathbb{P}\bigl(\Phi(X_N)>a\bigr) = \int
\mathbb{P}\bigl(\Phi(X_N^{n_0,x_0})>a\bigr) d\mu_0(x_0)}$ where $\mu_0$
denotes the law of $X_{n_0}$.
\end{rem} 

\section{General formulation of the Adaptive Multilevel Splitting algorithm}\label{sec:generalAMS}

\subsection{Context}\label{sec:context}

The goal is to estimate the  probability $p$ given
by~\eqref{eq:proba3}, in the regime where $p$ is small, which is for
example the case when $a$ is large.

It is convenient to introduce an auxiliary process $\bigl(Z_n\bigr)_{n_0\le n\le N}$, such that $Z_n=\bigl(n\Delta t,X_n^{n_0,x_0}\bigr)$. Indeed, let
\[
A=\left\{(T,x);~\Phi(x)\le a\right\}, B=\left\{(T,x);~\Phi(x)>a\right\},
\]
and define the associated stopping times
\[
\tau_A=\inf\left\{n\ge n_0,n\in\mathbb{N}_0;~Z_n\in A\right\}, \tau_B=\inf\left\{n\ge n_0,n\in\mathbb{N}_0;~Z_n\in B\right\}.
\]
Then the probability $p$ given by~\eqref{eq:proba3} can be rewritten as
\begin{equation}\label{eq:p_AMS}
p=\mathbb{P}\bigl(\Phi(X_N^{n_0,x_0})>a\bigr)=\mathbb{P}(\tau_B<\tau_A).
\end{equation}
We are then in position to build algorithms which fit in the
Generalized Adaptive Multilevel Splitting framework developed
in~\cite{BrehierGazeauGoudenegeLelievreRousset:16}, which ensure that
the obtained estimators of the probability~\eqref{eq:p_AMS} are unbiased.

For that, a score function, or reaction coordinate, $\xi$, needs to be given. Following the interpretation above, $\xi$ may depend on $z=(n\Delta t,x)$.

To run the algorithm and define simple unbiased estimators of $p$, only one requirement is imposed on the function $\xi$: there exists $\xi_{\rm max}$ such that
\[
B\subset \left\{z;~\xi(z)>\xi_{\rm max}\right\},
\]
which in the context of this article is rephrased as
\begin{equation}\label{eq:HYP}
\Phi(x)>a \implies \xi(T,x)>\xi_{\rm max}.
\end{equation}
The principle of splitting algorithm is then to write
$$\mathbb P(\tau_B < \tau_A)=\mathbb P(\tau_{\zeta_1}<\tau_A) \mathbb P(\tau_{\zeta_2}<\tau_A |\tau_{\zeta_1}<\tau_A) \mathbb P(\tau_{\zeta_3}<\tau_A |\tau_{\zeta_2}<\tau_A)
\ldots \mathbb P(\tau_B <\tau_A | \tau_{\xi_{\rm max}} <\tau_A) $$
for an increasing sequence of levels $(\zeta_q)_{q \ge 1}$, where
$\tau_\zeta=\inf\{n \ge n_0; \, \xi(Z_n)>\zeta\}$. If the levels are well
chosen, then the successive conditional probabilities $\mathbb P(\tau_{\zeta_{q+1}}<\tau_A |
\tau_{\zeta_q}<\tau_A)$ are easy to compute. The principle of the adaptive
multilevel splitting algorithm~\cite{CerouGuyader:07} is to choose the levels
adaptively, so that the successive conditional probabilities $\mathbb P(\tau_{\zeta_{q+1}}<\tau_A | \tau_{\zeta_q}<\tau_A)$ are constant and fixed. The levels constructed in the algorithm are then random.

\subsection{The Adaptive Multilevel Splitting algorithm}\label{sec:algo}

Before giving the detailed algorithm, let us roughly explain the main
steps (we also refer
to~\cite{BrehierGazeauGoudenegeLelievreRousset:16} for more details and intuition on the algorithm). In the initialization, one samples $n_{\rm rep}$ trajectories
following~\eqref{eq:SDE}-\eqref{eq:SDE-IC} and compute the score of
each trajectory, namely the maximum of $\xi$ attained along the
path. Then the algorithm proceeds as follows: one discards the
trajectory which has the smallest score and in order to keep the
number of trajectories fixed, a new one is created by choosing one of
the remaining trajectories at random, copying it up to the score of
the killed trajectory, and sampling the end of trajectory
independently from the past. This is called the partial resampling. One thus obtains a new ensemble of
$n_{\rm rep}$ trajectories on which one can iterate by again
discarding the the trajectory which has the smallest score. As the
iteration goes, one thus obtains trajectories with largest and largest
scores, and an estimate of the probability of interest is obtained as
$(1-1/n_{\rm rep})^{Q_{\rm iter}} P(\tau_B <\tau_A | \tau_{\xi_{\rm max}} <\tau_A)$ (notice that $(1-1/n_{\rm rep})$ is an estimate of the conditional probability
to reach level $\zeta_{q+1}$ conditionally to the fact that level $\zeta_q$
has been reached), where $Q_{\rm iter}$ is the number of iterations
required to reach the maximum level $\xi_{\rm max}$. In practice,
$P(\tau_B <\tau_A | \tau_{\xi_{\rm max}} <\tau_A)$ is estimated by the
proportion of trajectories which reach $B$ before $A$ at the last
iteration of the algorithm, namely when all the trajectories satisfy $\tau_{\xi_{\rm max}} <\tau_A$.

Actually, the algorithm has to be adapted in order to take into
account situations when more than one particle has the smallest
score, which happens with non zero probability for Markov chains. Let us now give the details of the AMS algorithm.

To simplify notation, in the sequel, the initial condition $x_0$ and
the time $n_0$ are omitted in the notation of the replicas.

\paragraph{Input}
\begin{itemize}
\item $n_{\rm rep}\in\mathbb{N}$, the number of replicas,
\item a score function $z=(n\Delta t,x)\mapsto\xi(n \Delta
  t,x)\in\mathbb{R}$ and a stopping level $\xi_{\rm max}\in\mathbb{R}$
  such that~\eqref{eq:HYP} is satisfied.
\end{itemize}

\paragraph{Initialization}
\begin{itemize}
\item Sample $\nrep$ independent realizations of the Markov process
\[
X^j=\bigl(X_m^j\bigr)_{n_0\le m\le N},~1\le j\le \nrep
\]
following the dynamics~\eqref{eq:SDE}--\eqref{eq:SDE-IC}.
\item Compute the score of each replica, $M^{j}=\underset{n_0\le m\le N}\max \xi(m\Delta t,X_m^j)$.
\item Compute the level ${\mathcal{Z}}=\underset{1\le j\le \nrep}\min M^j$.
\item Define $\mathcal{K}=\left\{j\in\left\{1,\ldots,\nrep\right\}~;~M^j=\mathcal{Z}\right\}$.
\item Set $q=0$, $\hat{p}=1$, $\mathcal{B}=1$.
\end{itemize}

\paragraph{Stopping criterion}

{\bf If} {$\mathcal{Z}\ge \xi_{\rm max}$ or ${\rm card}(\mathcal{K})=\nrep$}, {\bf then} set $\mathcal{B}=0$.

\paragraph{While $\mathcal{B}=1$}

\begin{itemize}
\item {\bf Update}
\begin{itemize}
\item $q \leftarrow q+1$ and $\hat{p} \leftarrow \hat{p}\cdot \bigl(1-\frac{{\rm card}(\mathcal{K})}{n}\bigr)$.
\end{itemize}
\item {\bf Splitting}
\begin{itemize}
\item Reindex the replicas, such that
\[
\begin{cases}
M^j=\mathcal{Z} \quad \text{if}~ j\in\left\{1,\ldots,{\rm card}(\mathcal{K})\right\}\\
M^j>\mathcal{Z} \quad \text{if}~j\in\left\{{\rm card}(\mathcal{K})+1,\ldots,\nrep\right\}.
\end{cases}
\]
\item For replicas with index $j\in\left\{1,\ldots,{\rm card}(\mathcal{K})\right\}$, sample labels $\ell_1,\ldots,\ell_{{\rm card}(\mathcal{K})}$, independently and uniformly in $\left\{{\rm card}(\mathcal{K})+1,\ldots,\nrep\right\}$.
\end{itemize}
\item {\bf Partial resampling}
\begin{itemize}
\item Remove the replicas with label $j\in\left\{1,\ldots,{\rm card}(\mathcal{K})\right\}$. 
\item For $j\in\left\{1,\ldots,{\rm card}(\mathcal{K})\right\}$, define $m_j=\inf\left\{m\in\left\{n_0,\ldots,N\right\}~;~\xi(X_m^{\ell_j})>\mathcal{Z}\right\}$.
\item For $m\in\left\{n_0,\ldots,m_j\right\}$, set $X_m^{j}=X_m^{\ell_j}$.
\item Sample a new trajectory $\bigl(X_m^j\bigr)_{m_j\le m\le N}$ with
  the Markov dynamics~\eqref{eq:SDE} driven by independent
  realizations of the Brownian motion.
\end{itemize}
\item {\bf Level computation}
\begin{itemize}
\item Compute the scores $M^{j}=\underset{n_0\le m\le N}\max \xi(m\Delta t,X_m^j)$.
\item Compute the level $\mathcal{Z}=\underset{1\le j\le \nrep}\min M^j$.
\item Define the set $\mathcal{K}=\left\{j\in\left\{1,\ldots,\nrep\right\}~;~M^j=\mathcal{Z}\right\}$.
\end{itemize}
\end{itemize}

\paragraph{Stopping criterion}

{\bf If} {$\mathcal{Z}\ge \xi_{\rm max}$ or ${\rm card}(\mathcal{K})=\nrep$}, {\bf then} set $\mathcal{B}=0$.

\paragraph{End while}
\paragraph{Update:} $\hat{p}\leftarrow \hat{p} \frac{1}{\nrep}\sum_{j=1}^{\nrep}\mathds{1}_{\Phi(X_N^j)\ge a}$.
\paragraph{Output:} $\hat{p}$ and $Q_{\rm iter}=q$.

\begin{rem}
We presented the algorithm in its simplest form. There are many
variants, see~\cite{BrehierGazeauGoudenegeLelievreRousset:16}. For
example, the killing level $\mathcal{Z}$ can be defined as $\mathcal{Z}=M^{(k)}$ where
$M^{(1)} \le M^{(2)} \le \ldots \le M^{(n_{\rm rep})}$ denotes an increasing relabelling of the scores $(M^j)_{1 \le j \le n_{\rm rep}}$ (order statistics).
\end{rem}

\subsection{Consistency result}

Let $\hat{p}$ and $Q_{\rm iter}$ be the outputs of the realization of the algorithm above. We quote the following result~\cite{BrehierGazeauGoudenegeLelievreRousset:16}, which states that the output $\hat{p}$ of the algorithm described in Section~\ref{sec:algo} above is an unbiased estimator of the probability given by~\eqref{eq:proba3}.
\begin{theo}\label{theo:unbias}
Let $\xi$ be a score function and $\xi_{\rm max}\in \mathbb{R}$ be
such that~\eqref{eq:HYP} is satisfied.
Let $\nrep\in\mathbb{N}$ be a given number of replicas. Assume that almost surely the algorithm stops after a finite number of iterations: $Q_{\rm iter}<\infty$ almost surely.

Then $\hat{p}$ is an unbiased estimator of the probability $p$ given by~\eqref{eq:p_AMS}:
\[
\mathbb{E}[\hat{p}]=\mathbb{P}\bigl(\Phi(X_N^{n_0,x_0})>a\bigr).
\]
\end{theo}

Note that if $\psi:\mathbb{R}^d\to\mathbb{R}$ is a function with
support included in $\left\{x;~\Phi(x)>a\right\}$, {\it i.e.}
$\psi(x)=0$ if $\Phi(x)\le a$, then an unbiased estimator of
$\mathbb{E}[\psi(X_N^{n_0,x_0})]$ is given by replacing the final
update in the algorithm above (see step $g.$), by
\[
\hat{p}\leftarrow \hat{p} \frac{1}{\nrep}\sum_{j=1}^{\nrep}\psi(X_N^j).
\]

The unbiasedness property is crucial in practice for the following two
reasons. First, it is very easy to parallelize the estimation of rare
events using this property. Indeed, since the estimator is unbiased
whatever the value of $n_{\rm rep}$, to get a convergent estimation,
one has simply to fix $n_{\rm rep}$ to a value which enables the
computation of $\hat{p}$ on a single CPU, and then to sample $M$
independent realizations of this $\hat{p}$, run in parallel. In the
large $M$ limit, one obtains a convergent estimator of the quantity of
interest by simply considering the average of the realizations of
$\hat{p}$. Second, the practical interest of Theorem~\ref{theo:unbias}
is that since $\mathbb E(\hat{p})$ is the same whatever the choice of the
numerical parameters (namely $n_{\rm rep}$ and $\xi$), one can compare
the results obtained with different choices to get confidence in the
result. For example, one can consider the confidence intervals
obtained with $M$ independent realizations of $\hat{p}$ for two
different choices of $\xi$, and check whether these confidence
intervals overlap or not.

\begin{rem}
In the algorithm described in Section~\ref{sec:algo}, the set
$\mathcal{K}$ defined in the initialization and in the level
computation steps may have a cardinal strictly larger than $1$, even
if the level $\mathcal{Z}$ is defined as the minimum of the scores over the
replicas. This simply means than more than one replica has a score
which is the smallest among the replicas. In the discrete-time
setting (namely for Markov chains), this happens with non zero
probability, and it requires an appropriate modification of the
original AMS algorithm, as described above,
see~\cite{BrehierGazeauGoudenegeLelievreRousset:16} for more details.

Notice that in particular, there is a possibility that the algorithm stops if ${\rm card}(\mathcal{K})=\nrep$, in which case there is an extinction of the system of replicas.
\end{rem}

\section{Choices of the score function}\label{sec:choice}

Let us now describe the various score functions we will consider in
order to estimate~\eqref{eq:proba3}.

\subsection{Vanilla score function and limitations}\label{sec:choice_vanilla}

The simplest choice consists in choosing the score function as given by
\begin{equation}\label{eq:vanillaSF}
\xi_{\rm std}(n\Delta t,x)=\Phi(x),
\end{equation}
with $\xi_{\rm max}=a$. In this case, the score function does not depend on the time variable.

If the conditional probability
\[
q=\mathbb{P}\bigl(\Phi(X_N^{n_0,x_0})>a~\big|~\underset{n_0\le n\le N}\max~\Phi(X_n^{n_0,x_0})>a\bigr)
\]
is small, the performance of this vanilla strategy may be
poor. Indeed, only a small proportion of the replicas have a non-zero
contribution to the value of the estimator of the
probability~\eqref{eq:proba3}. It may even happen that all the
replicas satisfy $\underset{n_0\le n\le N}\max~\Phi(X_n^{n_0,x_0})>a$
but that none of them satisfies $\Phi(X_N^{n_0,x_0})>a$. In that
situation, the algorithm returns $\hat{p}=0$ to estimate $p>0$. One of
the goals of this work is to construct score functions which
circumvent that issue: with the score functions introduced below, almost surely $\hat{p}\neq 0$.

\subsection{Time-dependent score functions}\label{sec:choice_new}

As discussed above, it is natural to design score functions $\xi$, which satisfy the following condition:
\begin{equation}\label{eq:condition}
\left\{\Phi(X_N^{n_0,x_0})>a\right\}=\left\{\underset{n_0\le n\le T}\max~\xi(n\Delta t,X_n^{n_0,x_0})>1\right\}.
\end{equation}
The choice of the value $\xi_{\rm max}=1$ on the right-hand side above
is arbitrary, but no generality is lost. Indeed, if a score function $\xi$ satisfying~\eqref{eq:condition} is used in the AMS algorithm above, with $\xi_{\rm max}=1$, at the last update, the ratio $\frac{1}{\nrep}\sum_{j=1}^{\nrep}\mathds{1}_{\Phi(X_N^j)\ge a}$ is identically equal to $1$, since all replicas satisfy $\underset{n_0\le n\le T}\max~\xi(n\Delta t,X_n^{n_0,x_0})>1$. In particular, by construction $\hat{p}\neq 0$ (provided $Q_{\rm iter}<\infty$).

As will be seen below, in practice it is more natural to identify functions $\tilde{\xi}$ taking values in $(-\infty,1]$, which satisfy the condition
\begin{equation}\label{eq:condition_bis}
\left\{\Phi(X_N^{n_0,x_0})>a\right\}=\left\{\underset{n_0\le n\le T}\max~\tilde{\xi}(n\Delta t,X_n^{n_0,x_0})=1\right\},
\end{equation}
instead of~\eqref{eq:condition}. To justify the use of the algorithm
in this case, observe that $\bar
\xi(t,x)=\tilde{\xi}(t,x)+\mathds{1}_{\tilde{\xi}(t,x)=1}$ then
satisfies~\eqref{eq:condition}. In addition, when running the
algorithm, choosing either $\bar \xi$ or $\tilde{\xi}$ exactly yields
the same result. We are thus in the setting where the unbiasedness
result Theorem~\ref{theo:unbias} applies. The score functions presented below will satisfy~\eqref{eq:condition_bis} instead of~\eqref{eq:condition}.

\begin{figure}
\includegraphics[scale=0.4]{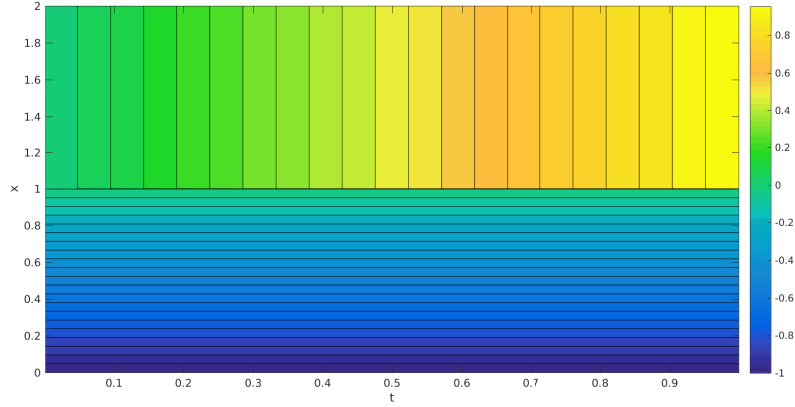}
\caption{Level lines of the score function $(t,x)\mapsto \xi_{\rm new}(t,x)$, with $\Phi(x)=x$, $a=1$, $T=1$.}
\label{fig:xi_new}
\end{figure}

One of the novelties of this article is the introduction of the following score function:
\begin{equation}\label{eq:scorefunction}
\xi_{\rm new}(n\Delta t,x)=\bigl(\Phi(x)-a\bigr)\mathds{1}_{\Phi(x)\le a}+\frac{n\Delta t}{N\Delta t}\mathds{1}_{\Phi(x)>a}.
\end{equation}
Observe that $\xi_{\rm new}$ takes values in $(-\infty,1]$, and that
$\xi_{\rm new}(n\Delta t,x)=1$ if and only if $n=N$ and $\Phi(x)\ge
a$. Thus the condition~\eqref{eq:condition_bis} is satisfied. We refer
to Figure~\ref{fig:xi_new} for a schematic representation of this score function.

Note that the score function defined by~\eqref{eq:scorefunction} only depends on the function $\Phi$, on the threshold $a$, and on the final time $T=N\Delta t$. It may thus be applied in any situation, but in some cases better score functions may be built upon using more information on the dynamics. The practical implementation is very simple.

Let us explain how the AMS algorithm proceeds when used with the score function~\eqref{eq:scorefunction}. Observe that
\[
\left\{\underset{n_0\le n\le N}\max~\xi_{\rm new}(n\Delta t,X_n^{n_0,x_0})\ge 0\right\}=\left\{\underset{n_0\le n\le N}\max~\Phi(X_n^{n_0,x_0})>a\right\}.
\]
The first iterations of the algorithm, up to reaching level $0$, are
thus devoted to construct $n_{\rm rep}$ replicas which satisfy the
weaker condition $\left\{\underset{n_0\le n\le
    N}\max~\Phi(X_n^{n_0,x_0})>a\right\}$. In other words, if the
stopping level $\xi_{\rm max}$ in the algorithm is set equal to $0$ instead of $1$, one thus recovers the vanilla AMS algorithm described above, applied with the score function $\xi(t,x)=\Phi(x)$ (independent of time $t$).

Compared with the vanilla score function, the AMS algorithm the new
score function does not stop when $\left\{\underset{n_0\le n\le
    N}\max~\Phi(X_n^{n_0,x_0})>a\right\}$. In terms of splitting, observe that this consists in writing
\[
\mathbb{P}\bigl(\Phi(X_N^{n_0,x_0})>a\bigr)=\mathbb{P}\bigl(\Phi(X_N^{n_0,x_0})>a~|~\underset{n_0\le n\le N}\max~\Phi(X_n^{n_0,x_0})>a\bigr) \mathbb{P}\bigl(\underset{n_0\le n\le N}\max~\Phi(X_n^{n_0,x_0})>a\bigr),
\]
and the remaining effort consists in estimating the conditional probability above.

More generally, observe that for every $n_1\in\left\{n_0,\ldots,N\right\}$,
\[
\left\{\underset{n_0\le n\le N}\max~\xi_{\rm new}(n\Delta t,X_n^{n_0,x_0})\ge \frac{n_1}{N}\right\}=\left\{\underset{n_1\le n\le N}\max~\Phi(X_n^{n_0,x_0})>a\right\}.
\]
The construction of the score function~\eqref{eq:scorefunction} is associated to the following family of nested events:
\begin{align*}
\left\{\underset{n_0\le n\le N}\max~\xi_{\rm new}(n\Delta t,X_n^{n_0,x_0})\ge 1\right\}&\subset \ldots \subset \left\{\underset{n_0\le n\le N}\max~\xi_{\rm new}(n\Delta t,X_n^{n_0,x_0})\ge \frac{n_1}{N}\right\}\\
&\subset  \ldots \subset \left\{\underset{n_0\le n\le N}\max~\xi_{\rm new}(n\Delta t,X_n^{n_0,x_0})\ge 0\right\},
\end{align*}
which equivalently may be rewritten as
\begin{align*}
\left\{\Phi(X_N^{n_0,x_0})>a\right\}&\subset \ldots \subset \left\{\underset{n_1\le n\le N}\max~\Phi(X_n^{n_0,x_0})>a\right\}\\
&\subset  \ldots \subset \left\{\underset{n_0\le n\le N}\max~\Phi(X_n^{n_0,x_0})>a\right\}.
\end{align*}
In the selection procedure, the intervals $[n_1\Delta t,N\Delta t]$ are iteratively reduced (by increasing the left end point of the interval), until they ultimately contain only the point $N\Delta t$ (the right end point of the interval which remains fixed).

\begin{figure}
\includegraphics[scale=0.4]{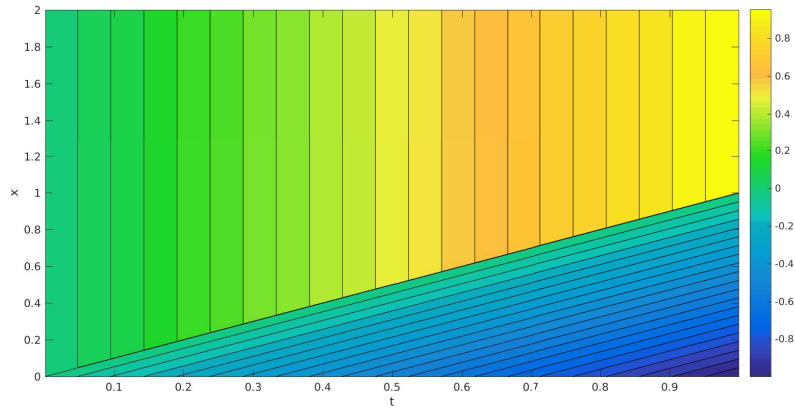}
\caption{Level lines of the score function $(t,x)\mapsto \xi_{\rm
    new}^{\mathbf{a}}(t,x)$, with $\mathbf{a}(t)=atT^{-1}$, and
  $\Phi(x)=x$, $a=1$, $T=1$.}
\label{fig:xi_new_a}
\end{figure}

To conclude, we mention that the construction given by~\eqref{eq:scorefunction} can be generalized as follows. Let $\mathbf{a}:[0,T]\to\mathbb{R}$ be a non-decreasing function, such that $\mathbf{a}(T)=a$. Define
\begin{equation}\label{eq:xinewa}
\xi_{\rm new}^{\mathbf{a}}(n\Delta t,x)=\bigl(\Phi(x)-\mathbf{a}(n\Delta t)\bigr)\mathds{1}_{\Phi(x)\le \mathbf{a}(n\Delta t)}+\frac{n\Delta t}{N\Delta t}\mathds{1}_{\Phi(x)>\mathbf{a}(n\Delta t)}.
\end{equation}
The score function $\xi_{\rm new}$ defined by~\eqref{eq:scorefunction}
is a particular case of~\eqref{eq:xinewa}, with
$\mathbf{a}(t)=a$. Optimizing the choice of the function $\mathbf{a}$
may help improve the efficiency of the
algorithm. Notice that condition~\eqref{eq:condition_bis} is satisfied with $\xi=\xi_{\rm new}^{\mathbf{a}}$.  We refer to Figure~\ref{fig:xi_new_a} for a schematic representation of this
score function.

\subsection{The optimal score function: the committor function}\label{sec:committor}

For the general setting presented in Section~\ref{sec:context} where
ones want to estimate $\mathbb P(\tau_B < \tau_A)$, the committor function is
defined as $z \mapsto \mathbb P^z(\tau_B <
\tau_A)$, where the upperscript $z$ refers to the initial
condition of the process $Z$. In~\cite{CerouDelyonGuyaderRousset:18}, it is shown that, in
a continous-time setting, the asymptotic variance (as the number of
replicas $\nrep$ goes to infinity) of AMS algorithm is minimized when
using the committor function as the score function. It is thus
interesting to look at what the committor function looks like in our context.

In our context, the committor function is given by
\begin{equation}\label{eq:com}
\xi_{\rm com}(n\Delta t,x)=\mathbb{P}\bigl(\Phi(X_N^{n,x})>a\bigr).
\end{equation}
For the discussion, it is more convenient to consider the
continuous-time version $\xi_{\rm  com}(t,x)=\mathbb{P}\bigl(\Phi(X_T^{t,x})>a\bigr)$,  for
$t\in[0,T]$, that we still denote $\xi_{\rm com}$ with a slight abuse of notation. Recall that
$u_a(t,x)=\mathbb{P}\bigl(\Phi(X_T^{t,x})>a\bigr) = \xi_{\rm  com}(t,x)$ satisfies the
Kolmogorov backward equation~\eqref{eq:KB}, as explained in Section~\ref{sec:setting}.

As mentioned above, the asymptotic variance (as the number of replicas
$\nrep$ goes to infinity) of AMS algorithm is minimized when using the
committor function as the score
function~\cite{CerouDelyonGuyaderRousset:18}. The
asymptotic variance is then $\frac{-p^2\log(p)}{\nrep}$, where $p$ is
the probability which is estimated.  The analysis of the AMS algorithm
in the ideal case, {\it i.e.} when using the committor function as a
the score function, has been performed in many
works\cite{BrehierLelievreRousset:15,BrehierGoudenegeTudela:16,Brehier:16,GuyaderHentgartnerMatzner-Lober:11}.

Of course, in practice, the committor function is unknown and the
asymptotic variance depends on the chosen score function. It has been
proved~\cite{CerouDelyonGuyaderRousset:18} that the asymptotic
variance is always bounded from above by $\frac{2p(1-p)}{\nrep}$, for
any choice of the score function, where we recall that the asymptotic
variance of the vanilla Monte-Carlo method is
$\frac{p(1-p)}{\nrep}$. This can be seen as a sign of the robustness of the AMS approach to estimate rare event probability (contrary to importance sampling method which may result in a dramatic increase of the asymptotic variance compared with the vanilla Monte-Carlo method).

\begin{figure}
\includegraphics[scale=0.4]{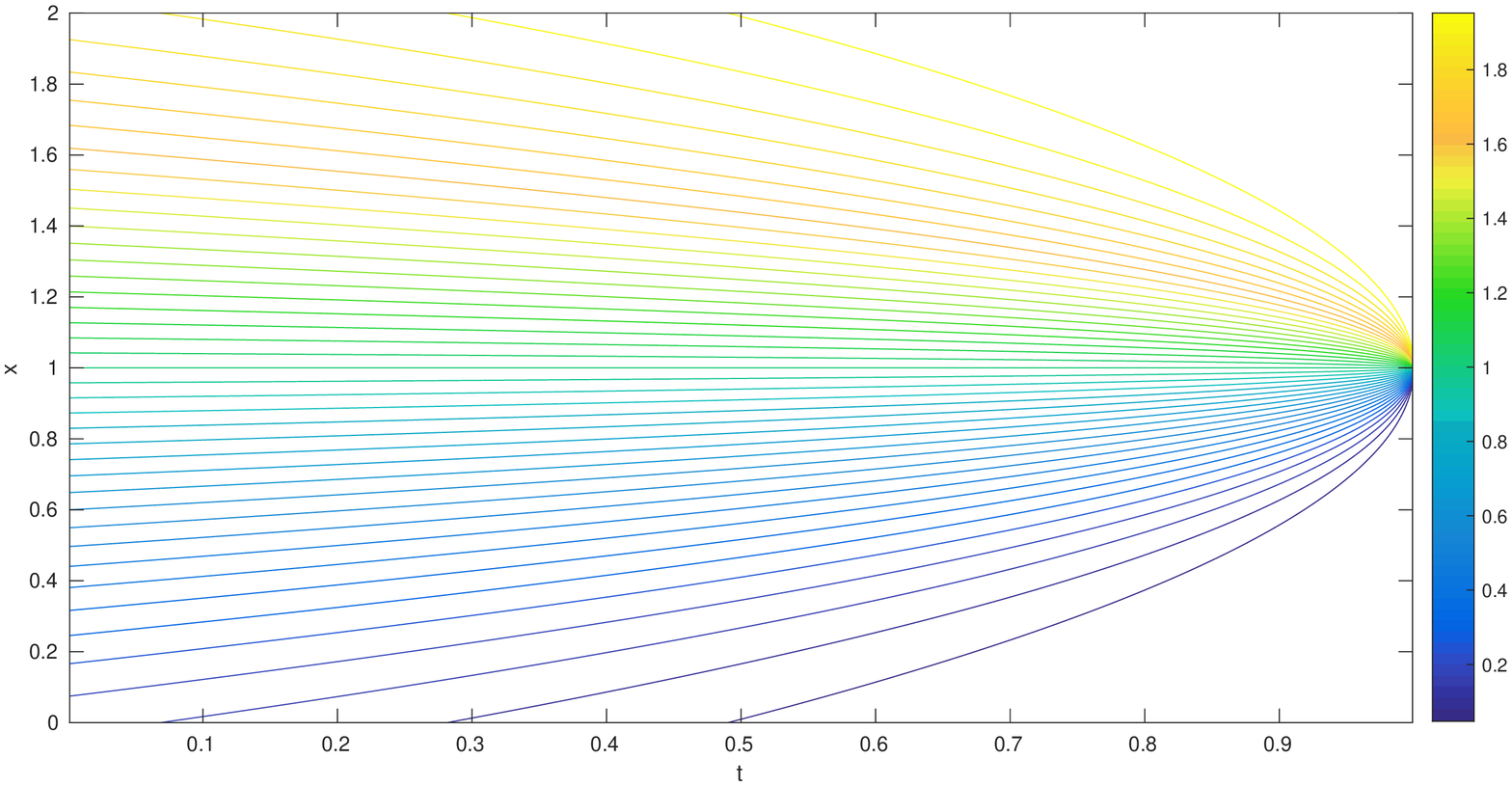}
\caption{Level lines of the committor function $(t,x)\mapsto \xi_{\rm com}(t,x)$, in the Brownian Motion case $X(t)=B(t)$, with $a=1$ and $T=1$.}
\label{fig:committor}
\end{figure}

For simple Gaussian models, namely when $X$ is a Brownian Motion, an
Ornstein-Uhlenkeck process, or a drifted Brownian Motion, it is
possible to compute analytically the committor function. This is useful to validate algorithms on test cases, as will be illustrated in Section~\ref{sec:numerical}. Figure~\ref{fig:committor} represents the level lines of the committor function for a one-dimensional Brownian Motion (with $T=1$ and $a=1$). In that case,
\[
\xi_{\rm com}(t,x)=1-F\left(\frac{a-x}{\sqrt{T-t}}\right),
\]
where $F$ is the cumulative distribution function of the standard
Gaussian distribution $\mathcal{N}(0,1)$, to be compared to the level
sets of $\xi_{\rm new}$ and $\xi^{\mathbf a}_{\rm new}$ on Figures~\ref{fig:xi_new}
and~\ref{fig:xi_new_a}. This form leads to define other families of appropriate score functions:
\[
\xi(t,x)=1-F\bigl(\phi(t,x)\bigr)
\]
where $\phi(t,x)\underset{t\to\infty}\to
(-\infty)\mathds{1}_{\Phi(x)>a}+(+\infty)\mathds{1}_{\Phi(x)< a}$. But
the efficiency depends a lot on the choice of $\phi$. In practice, we
did not observe much gain in our numerical experiments, compared to the score function $\xi_{\rm
  new}$ introduced in the previous section.

Let us mention that various techniques have been proposed in order to
approximate the committor function, in particular in the context of
importance sampling techniques for rare events, since the committor
function also gives the optimal change of measure. If diffusions with vanishing noise are
considered~\cite{Vanden-EinjdenWeare2012,dupuis-spiliopoulos-wang-12,dupuis-spiliopoulos-zhou-15}, solutions of associated
Hamilton-Jacobi equations are good candidates to estimate the
committor function. See also~\cite{schuette-winkelman-hartmann-12} for approximations based on
coarse-grained models. Whether such constructions are possible when
considering temporal averages, instead of the terminal value of the
process is unclear.

\section{Numerical simulations}\label{sec:numerical}

Let $p$ denote the rare event probability of interest. An estimator of
$p$ is calculated as the empirical average over independent
realizations of the AMS algorithm, given a choice of score function
$\xi$. The main objective of this section is to investigate the
behavior of the algorithm when choosing $\xi=\xi_{\rm new}$ given
by~\eqref{eq:scorefunction}. A comparison with the vanilla score
function $\xi=\xi_{\rm std}$ given by~\eqref{eq:vanillaSF} is provided.

Let $\nrep\in\mathbb{N}$ denote the number of replicas, and let $M\in\mathbb{N}$, and $\bigl(\hat{p}_m\bigr)_{1\le m\le M}$ be the output probabilities of $M$ independent realizations of the AMS algorithm. We report the values of the empirical average
\[
\hat{p}=\frac{1}{M}\sum_{m=1}^{M}\hat{p}_m,
\]
and of the empirical variance $\hat{\sigma}^2=\frac{1}{M-1}\sum_{m=1}^{M}(\hat{p}_m-\hat{p})^2$. Confidence intervals are computed as
\[
\left[\hat{p}-\frac{1.96 \hat{\sigma}}{\sqrt{M}},\hat{p}+\frac{1.96\hat{\sigma}}{\sqrt{M}}\right],
\]
assuming that the number of realizations $M$ is sufficiently large to use the Gaussian, Central Limit Theorem, regime.

Recall that $\mathbb{E}[\hat{p}]=\mathbb{E}[\hat{p}_1]=p$, whatever
the choice of the score function~$\xi$ and of the number of
replicas~$\nrep$, thanks to Theorem~\ref{theo:unbias}. The variance of the estimator and thus the efficiency strongly depends on $\xi$. In the experiments below, the
empirical variance $\hat{\sigma}^2$ is compared with the optimal
asymptotic variance $\frac{-p^2\log(p)}{n_{\rm rep}}$ for (adaptive)
multilevel splitting algorithms, which is obtained in the regime
$\nrep\to \infty$, when choosing the (unknown in general) committor
function $\xi=\xi_{\rm com}$ as the score function. The difference
between the empirical variance and the optimal one can be seen as a
measure of how far the chosen score function is from the committor.

In some of the numerical experiments below, the conditional probability
\begin{equation}\label{eq:q}
q=\mathbb{P}\left(\Phi(X_N)>a~\big|~\underset{0\le n\le N}\max~\Phi(X_n)>a\right),
\end{equation}
is also estimated by
\[
\hat{q}=\frac{\hat{p}}{\hat{p}_{\rm max}},
\]
where $\hat{p}_{\rm max}=\frac{1}{M}\sum_{m=1}^{M}\hat{p}_{\rm max,m}$ is the estimator of the probability
\begin{equation}\label{eq:pmax}
p_{\rm max}=\mathbb{P}\left(\underset{0\le n\le N}\max~\Phi(X_n)>a\right),
\end{equation}
which is estimated  using the vanilla score function $\xi=\xi_{\rm
  std}$.

\subsection{Validation using two Gaussian models}

In this section, we validate the AMS alogorithm with various score
functions on simple models for which the probability of the rare event
is known with arbitrary precision.

\subsubsection{Brownian Motion}

We follow here numerical experiments
from~\cite{Vanden-EinjdenWeare2012}. Let $d=1$, and consider the diffusion process given by
\[
dX(t)=\sqrt{2\beta^{-1}}dW(t),~X(0)=0.1,
\]
where $\bigl(W(t)\bigr)_{t\ge 0}$ is a standard real-valued Wiener process.

The dynamics is discretized using the explicit Euler-Maruyama method,
with time-step size $\Delta t=10^{-3}$ (notice that the numerical
scheme gives here the exact solution):
\[
X_{n+1}=X_n+\sqrt{2\beta^{-1}\Delta t}\zeta_n,
\]
with $X_0=X(0)=0.1$, where $\big(\zeta_{n}\bigr)_{0\le n\le N}$ are independent standard Gaussian random variables.

The goal is to estimate the probability
\[
p=\mathbb{P}\bigl(|X_N|>1\bigr).
\]
This corresponds with the choice $\Phi(x)=|x|$, $a=1$, $T=1$ so that
$N=T/\Delta t= 10^3$.

Since $X(t)=X(0)+W(t)$, the law of $X(t)$ is a Gaussian distribution, and the value of $\mathbb{P}(|X(1)|\ge 1)$ can be computed exactly in terms of the cumulative distribution function of the standard Gaussian distribution.

\bigskip

Two numerical experiments are reported below, using $\xi=\xi_{\rm new}$. First, in Table~\ref{tab:BM1}, the number of replicas is set equal to $n_{\rm rep}=10^2$, and the empirical average is computed over $M=10^4$ independent realizations of the algorithm. Second, in Table~\ref{tab:BM2}, the number of replicas is set equal to $n_{\rm rep}=10^3$, and the empirical average is computed over $M=10^3$ independent realizations of the algorithm.

\begin{table}[htbp]
\begin{center}
\begin{tabular}{||c||c|c||c||c|c||}
\hline
$\beta$ & $\hat{p}$ & $p$ ($\Delta t=0)$ & confidence interval & $\hat{\sigma}^2$ & $\frac{-p^2\log(p)}{n_{\rm rep}}$ \\
\hline
$2$ & $3.199~10^{-1}$ & $3.197~10^{-1}$ & $[3.192~10^{-1},3.205~10^{-1}]$ & $1.256~10^{-3}$ & $1.166~10^{-3}$ \\
$4$ & $1.613~10^{-1}$ & $1.614~10^{-1}$ & $[1.608~10^{-1},1.617~10^{-1}]$ & $5.275~10^{-4}$ & $4.751~10^{-4}$  \\
$8$ & $4.978~10^{-2}$ & $4.983~10^{-2}$ & $[4.958~10^{-2},4.998~10^{-2}]$ & $1.030~10^{-4}$ & $7.447~10^{-5}$  \\
$16$ & $6.395~10^{-3}$ & $6.386~10^{-3}$ & $[6.353~10^{-3},6.436~10^{-3}]$ & $4.449~10^{-6}$ & $2.061~10^{-6}$  \\
$32$ & $1.631~10^{-4}$ & $1.645~10^{-4}$ & $[1.611~10^{-4},1.651~10^{-4}]$ & $1.052~10^{-8}$ & $2.358~10^{-9}$  \\
$64$ & $1.787~10^{-7}$ & $1.782~10^{-7}$ & $[1.721~10^{-7},1.854~10^{-7}]$ & $1.139~10^{-13}$ & $4.935~10^{-15}$  \\
$128$ & $2.501~10^{-13}$ & $3.011~10^{-13}$ & $[1.766~10^{-13},3.235~10^{-13}]$ & $1.403~10^{-23}$ & $2.614~10^{-26}$\\
\hline
\end{tabular}
\caption{Brownian Motion, $n_{\rm rep}=10^2$, $M=10^4$.}
\label{tab:BM1}
\end{center}
\end{table}

\begin{table}[htbp]
\begin{center}
\begin{tabular}{||c||c|c||c||c|c||}
\hline
$\beta$ & $\hat{p}$ & $p$ ($\Delta t=0)$ & confidence interval & $\hat{\sigma}^2$ & $\frac{-p^2\log(p)}{n_{\rm rep}}$ \\
\hline
$2$ & $3.196~10^{-1}$ & $3.197~10^{-1}$ & $[3.188~10^{-1},3.203~10^{-1}]$ & $1.362~10^{-4}$ & $1.166~10^{-4}$ \\
$4$ & $1.617~10^{-1}$ & $1.614~10^{-1}$ & $[1.612~10^{-1},1.621~10^{-1}]$ & $5.456~10^{-5}$ & $4.751~10^{-5}$  \\
$8$ & $4.983~10^{-2}$ & $4.983~10^{-2}$ & $[4.963~10^{-2},5.002~10^{-2}]$  & $9.815~10^{-6}$  & $7.447~10^{-6}$  \\
$16$ & $6.411~10^{-3}$ & $6.386~10^{-3}$ & $[6.371~10^{-3},6.451~10^{-3}]$ & $4.242~10^{-7}$ & $2.061~10^{-7}$  \\
$32$ & $1.634~10^{-4}$ & $1.645~10^{-4}$ & $[1.614~10^{-4},1.655~10^{-4}]$ & $1.063~10^{-9}$ & $2.358~10^{-10}$  \\
$64$ & $1.800~10^{-7}$ & $1.782~10^{-7}$ & $[1.740~10^{-7},1.860~10^{-7}]$ & $9.360~10^{-15}$ & $4.935~10^{-16}$  \\
$128$ & $3.045~10^{-13}$ & $3.011~10^{-13}$ & $[2.426~10^{-13},3.664~10^{-13}]$  & $9.986~10^{-25}$ & $2.614~10^{-27}$  \\
\hline
\end{tabular}
\caption{Brownian Motion, $n_{\rm rep}=10^3$, $M=10^3$.}
\label{tab:BM2}
\end{center}
\end{table}

These numerical experiments validate the algorithm using $\xi=\xi_{\rm new}$ in the case of a one-dimensional Brownian Motion. The empirical variance $\hat{\sigma}^2$ is much smaller than $\frac{p(1-p)}{\nrep}$ which would be obtained using a naive Monte-Carlo strategy (using $\nrep$ independent replicas). It is observed that the ratio between the empirical and the optimal variances increase when $p$ decreases, but in practice this increase only has a limited impact and the new algorithm remains effective.

\subsubsection{Ornstein-Uhlenbeck}

We consider here an example taken from~\cite{WoutersBouchet:16}.
Let $d=1$, and consider the diffusion process given by
\[
dX(t)=-X(t)dt+dW(t),~X(0)=0.
\]

The dynamics is discretized using the explicit Euler-Maruyama method, with time-step size $\Delta t=10^{-3}$:
\[
X_{n+1}=X_n-\Delta tX_n+\sqrt{\Delta t}\zeta_n,
\]
with $X_0=X(0)=0$, where $\big(\zeta_{n}\bigr)_{0\le n\le N}$ are independent standard Gaussian random variables.

The goal is to estimate the probability
\[
p=\mathbb{P}\bigl(X_N>a\bigr),
\]
for different values of $a$. This corresponds with the choice
$\Phi(x)=x$. The value of $T$ is set to $T=2$ so that $N=T/\Delta t=2000$.

In this numerical experiment, the number of replicas is set equal to $n_{\rm rep}=10^2$, and the empirical average is computed over $M=10^4$ independent realizations of the algorithm.
The estimator $\hat{p}_{\rm new}$ of $p$ and the empirical variance
$\hat{\sigma}_{\rm new}^2$ are obtained using the score function
$\xi=\xi_{\rm new}$. The estimator $\hat{p}_{\rm std}$ and the empirical variance $\hat{\sigma}_{\rm std}^2$ are obtained using the vanilla splitting strategy, with reaction coordinate $\xi=\xi_{\rm std}$.
The value of the probability $p$ for the continuous time process, and the optimal variance $\frac{-p^2\log(p)}{n_{\rm rep}}$ are also reported for comparison.

\begin{table}[htbp]
\begin{center}
\begin{tabular}{||c||c|c|c||c|c|c||}
\hline
$a$ & $\hat{p}_{\rm new}$ & $\hat{p}_{\rm std}$  & $p$ ($\Delta t=0$) & $\hat{\sigma}_{\rm new}^2$ & $\hat{\sigma}_{\rm std}^2$ & $\frac{-p^2\log(p)}{n_{\rm rep}}$ \\
\hline
$2.8$ & $3.216~10^{-5}$ & $3.252~10^{-5}$ & $3.213~10^{-5}$ & $2.377~10^{-10}$  & $3.327~10^{-10}$ & $1.068~10^{-10}$\\
$2.9$ & $1.756~10^{-5}$ & $1.728~10^{-5}$ & $1.742~10^{-5}$& $7.917~10^{-10}$ & $1.009~10^{-10}$ & $3.325~10^{-11}$\\
$3.0$ & $9.341~10^{-6}$ & $9.300~10^{-6}$ & $9.260~10^{-6}$ & $2.411~10^{-11}$ & $5.832~10^{-11}$ & $9.937~10^{-12}$\\
$3.1$ & $4.857~10^{-6}$ & $4.826~10^{-6}$ & $4.827~10^{-6}$ & $7.011~10^{-12}$ & $8.482~10^{-12}$ & $2.852~10^{-12}$\\
$3.2$ & $2.486~10^{-6}$ & $2.449~10^{-6}$ & $2.468~10^{-6}$ & $1.984~10^{-12}$ & $2.475~10^{-12}$ & $7.864~10^{-13}$\\
\hline
\end{tabular}
\caption{Ornstein-Uhlenbeck, $T=2$. Comparison of the new and of the vanilla splitting algorithms. $n_{\rm rep}=10^2$, $M=10^4$.}
\label{tab:0U2}
\end{center}
\end{table}

The results of numerical experiments with $T=4$ and $T=8$ are reported below. The quantity
\[
\hat{r}=\frac{1}{M}\sum_{m=1}^{M}\mathds{1}_{\hat{p}_m>0}
\]
is also reported, when the vanilla score function is used. This is the
proportion of the independent realizations of the algorithm which
contribute in the empirical average. This proportion depends on the
conditional probability $q$ (see~\eqref{eq:q}): it may happen that the $n_{\rm rep}$
replicas obtained at the final iteration all satisfy $X_N\le 1$, even
if by construction they all satisfy $\underset{0\le n\le
  N}\max~X_n>1$. However, by construction (except if extinction of the
system of replicas happens, which has not been observed in this
experiment), if the new score function is used, $\hat{r}$ is
identically equal to $1$. Observe that, when $T$ goes to infinity, by
ergodicity of the process, $p_{\rm max}\to 1$ (see~\eqref{eq:pmax} for
the definition of $p_{\rm max}$), whereas $q$ and $p$ converge to a
non trivial probability. Thus, when $T$ goes to infinity, it is expected that $\hat{r}$ will be equal to $0$ if $M$ is too small, when using the vanilla strategy with $\xi=\xi_{\rm std}$. The conditional probability $q$ is also estimated: $\hat{q}=0.07$ when $T=2$, $\hat{q}=0.02$ when $T=4$, $\hat{q}=0.01$ when $T=8$.

\begin{table}[htbp]
\begin{center}
\begin{tabular}{||c||c|c|c||c|c|c||c||}
\hline
$a$ & $\hat{p}_{\rm new}$ & $\hat{p}_{\rm std}$  & $p$ ($\Delta t=0$) & $\hat{\sigma}_{\rm new}^2$ & $\hat{\sigma}_{\rm std}^2$ & $\frac{-p^2\log(p)}{n_{\rm rep}}$ & $\hat{r}$ \\
\hline
$2.8$ & $3.743~10^{-5}$  & $3.789~10^{-5}$ & $3.740~10^{-5}$ &$7.106~10^{-10}$ & $1.089~10^{-9}$  & $1.426~10^{-10}$ & $0.85$ \\
$2.9$ & $2.052~10^{-5}$ & $2.069~10^{-5}$ & $2.049~10^{-5}$ & $2.438~10^{-10}$ & $3.484~10^{-10}$  & $4.532~10^{-11}$ & $0.83$\\
$3.0$ & $1.104~10^{-5}$ & $1.124~10^{-5}$ & $1.101~10^{-5}$ &  $7.540~10^{-11}$ & $1.147~10^{-10}$ & $1.384~10^{-11}$ & $0.81$\\
$3.1$ & $5.822~10^{-6}$ & $5.856~10^{-6}$ & $5.805~10^{-6}$ & $2.412~10^{-11}$ & $3.584~10^{-11}$ & $4.062~10^{-12}$ & $0.78$\\
$3.2$ & $3.022~10^{-6}$ & $3.016~10^{-6}$ & $3.002~10^{-6}$ & $7.452~10^{-12}$  & $9.964~10^{-12}$ & $1.146~10^{-12}$ & $0.75$\\
\hline
\end{tabular}
\caption{Ornstein-Uhlenbeck, $T=4$. Comparison of the new and of the vanilla splitting algorithms. $n_{\rm rep}=10^2$, $M=10^4$.}
\label{tab:OU4}
\end{center}
\end{table}

\begin{table}[htbp]
\begin{center}
\begin{tabular}{||c||c|c|c||c|c|c||c||}
\hline
$a$ & $\hat{p}_{\rm new}$ & $\hat{p}_{\rm std}$  & $p$ ($\Delta t=0$) & $\hat{\sigma}_{\rm new}^2$ & $\hat{\sigma}_{\rm std}^2$ & $\frac{-p^2\log(p)}{n_{\rm rep}}$ & $\hat{r}$ \\
\hline
$2.8$ & $3.792~10^{-5}$ & $3.714~10^{-5}$ & $3.751~10^{-5}$ & $1.696~10^{-9}$ & $2.435~10^{-9}$ & $1.434~10^{-10}$ & $0.54$\\
$2.9$ & $2.036~10^{-5}$ & $2.071~10^{-5}$ & $2.055~10^{-5}$ & $5.439~10^{-10}$ & $8.289~10^{-10}$  & $4.557~10^{-11}$ & $0.51$\\
$3.0$ & $1.103~10^{-5}$ & $1.128~10^{-5}$ & $1.104~10^{-5}$ & $1.843~10^{-10}$  & $2.745~10^{-10}$ & $1.392~10^{-11}$ & $0.49$\\
$3.1$ & $5.915~10^{-6}$ & $5.968~10^{-6}$ & $5.824~10^{-6}$ & $5.599~10^{-11}$ & $8.457~10^{-11}$ & $4.089~10^{-12}$ & $0.47$\\
$3.2$ & $2.978~10^{-6}$ & $3.022~10^{-6}$ & $3.013~10^{-6}$ & $1.639~10^{-11}$ & $2.368~10^{-11}$ & $1.154~10^{-12}$ & $0.43$\\
\hline
\end{tabular}
\caption{Ornstein-Uhlenbeck, $T=8$. Comparison of the new and of the vanilla splitting algorithms. $n_{\rm rep}=10^2$, $M=10^4$.}
\label{tab:OU8}
\end{center}
\end{table}

To conclude this section, note that on this example, the AMS
algorithms applied with the vanilla and
the new score functions have a similar quantitative behavior in terms
of asymptotic variance. However, their qualitative properties are
different. When the conditional probability $q$ gets small, the
advantage of the new score function is the fact that the output
$\hat{p}$ is always positive, so that even with a few realizations,
one gets a rough but informative approximation of the target probability.

\subsection{Drifted Brownian Motion}

We here considers numerical examples taken from~~\cite{BrehierGazeauGoudenegeLelievreRousset:arxiv,RollandSimonnet:15}. Let $d=1$, and consider the diffusion process given by
\[
dX_t=-\alpha dt+\sqrt{2\beta^{-1}}dW(t),~X(0)=0.
\]
The dynamics is discretized using the explicit Euler-Maruyama method,
with time-step size $\Delta t=10^{-2}$ (which gives again the exact
solution in this simple case):
\[
X_{n+1}=X_n-\alpha\Delta t+\sqrt{2\beta^{-1}\Delta t}\zeta_n,
\]
with $X_0=X(0)=0$, where $\big(\zeta_{n}\bigr)_{0\le n\le N}$ are independent standard Gaussian random variables.

The goal is to estimate the probability
\[
p=\mathbb{P}\bigl(X_N>1\bigr),
\]
thus $\Phi(x)=x$, $a=1$. One considers the final time $T=1$, so that
$N=T/\Delta t=10^2$. The value of $\alpha$ is set equal to
$\alpha=4$. As in the previous example, the value of $p$ is easy to
get using the fact that $X_N$ is Gaussian.

In this numerical experiment, see Table~\ref{tab:dBM}, three choices of score functions are considered. The number of replicas is $\nrep=10^3$. First, the estimator $\hat{p}_{\rm new}$ and the empirical variance $\hat{\sigma}_{\rm new}^2$ are obtained using $\xi=\xi_{\rm new}$, with a sample size $M=4.10^4$. Second, the estimator $\hat{p}_{\rm new,\mathbf{a}}$ and the empirical variance $\hat{\sigma}_{\rm new,\mathbf{a}}^2$ are obtained using $\xi=\xi_{\rm new}^{\mathbf{a}}$ with $\mathbf{a}(t)=\frac{at}{T}$, with a sample size $M=4.10^5$. Finally, the estimator $\hat{p}_{\rm std}$ and the empirical variance $\hat{\sigma}_{\rm std}^2$ are obtained using $\xi=\xi_{\rm std}$, with a sample size $M=4.10^4$. The sample sizes are chosen such that the total computational cost is of the same order for the three methods.

{\footnotesize
\begin{table}[htbp]
\begin{tabular}{||c||c|c|c|c||c|c|c|c||c||}
\hline
$\beta$ & $\hat{p}_{\rm new}$ & $\hat{p}_{\rm new,\mathbf{a}}$ & $\hat{p}_{\rm std}$ & $p$ & $\hat{\sigma}_{\rm new}^2$ & $\hat{\sigma}_{\rm new,\mathbf{a}}^2$ & $\hat{\sigma}_{\rm std}^2$ & $\frac{-p^2\log(p)}{n_{\rm rep}}$ & $\hat{r}$  \\
\hline
$1$ & $2.037~10^{-4}$ & $2.036~10^{-4}$ & $2.033~10^{-4}$ & $2.035~10^{-4}$ & $1.572~10^{-9}$ & $5.734~10^{-9}$ & $3.525~10^{-9}$ & $3.519~10^{-10}$ & $0.99$\\
$2$ & $2.843~10^{-7}$ & $2.870~10^{-7}$ & $2.878~10^{-7}$ & $2.867~10^{-7}$ & $4.714~10^{-14}$ & $2.348~10^{-12}$ & $9.527~10^{-14}$ & $1.238~10^{-15}$ & $0.69$ \\
$3$ & $4.613~10^{-10}$ & $4.325~10^{-10}$ & $4.705~10^{-10}$ & $4.571~10^{-10}$ & $1.817~10^{-18}$ & $2.197~10^{-17}$ & $3.084~10^{-18}$ & $4.493~10^{-21}$ & $0.11$ \\
$4$ & $7.620~10^{-13}$ & $6.975~10^{-13}$ & $7.582~10^{-13}$ & $7.687~10^{-13}$ & $5.034~10^{-23}$ & $4.388~10^{-22}$ & $8.343~10^{-23}$ & $1.648~10^{-26}$ & $0.01$\\
\hline
\end{tabular}
\caption{Drifted Brownian Motion, $\beta\in\left\{1,2,3,4\right\}$, Comparison of two versions of the new splitting algorithm and of the vanilla splitting algorithm.}
\label{tab:dBM}
\end{table}
}
Since the sample size is not the same for the three examples of score
functions in Table~\ref{tab:dBM}, the values of the empirical
variances $\hat{\sigma}^2$ should be taken with care when comparing
the methods. One would rather compare the values of
$\frac{\hat{\sigma}^2}{M}$. Then one concludes that the best
performance is obtained when using $\xi=\xi^{\mathbf a}_{\rm new}$. The
vanilla strategy, with $\xi=\xi_{\rm std}$, seems to behave
quantitatively the same as when $\xi=\xi_{\rm new}$. However, the
values of the proportion $\hat{r}$ of realizations such that
$\hat{p}_m\neq 0$ is not zero is also reported, when $\xi=\xi_{\rm
  std}$ (by construction, $\hat{r}=1$ for the first two cases). This
means that if $M$ was decreased (for instance, $M$ of the order $10^2$
for $\beta=4$), then the output of the experiment would be $\hat{p}_{\rm std}=0$.

As a consequence, the new algorithm clearly overcomes the limitation
of the vanilla strategy when $\xi=\xi_{\rm std}$. However, the score functions are far from being optimal, as revealed by the comparison with the optimal variance.

\subsection{Temporal averages for an Ornstein-Uhlenbeck process}

In this section, we consider an example taken from~\cite{FerreTouchette:18}.
Consider the one-dimensional Ornstein-Uhlenbeck process $X$, which is the solution of the SDE
\[
dX(t)=-X(t)dt+\sqrt{2\beta^{-1}}dW(t),~X(0)=0,
\]
and define the temporal average
\[
Y(T)=\frac{1}{T}\int_{0}^{T}X(s)ds.
\]
More generally, set $Y(t)=\frac{1}{t}\int_{0}^{t}X(s)ds$, for $t\in(0,T]$ and $Y(0)=0$.

The discretization is performed using the explicit Euler-Maruyama method, with time-step size $\Delta t=5.~10^{-3}$: for $n\in\left\{0,\ldots,N\right\}$ with $N\Delta t=T$,
\[
X_{n+1}=(1-\Delta t)X_n+\sqrt{2\beta^{-1}\Delta t}\zeta_n,\quad Y_{n}=\frac{1}{n}\sum_{m=1}^{n}X_m,
\]
with $X_0=Y_0=0$. Note that $Y$ satisfies a recursion formula $Y_{n+1}=(1-\frac{1}{n+1})Y_n+\frac{1}{n+1}X_{n+1}$.

The number of replicas is set equal to $\nrep=10^3$ and the sample size to compute empirical averages is $M=10^2$.

In this section, the probability which is estimated is
\[
p(T,a)=\mathbb{P}(Y_N>a).
\]
 The associated estimator is denoted by $\hat{p}(T,a)$ and the empirical variance by $\hat{\sigma}^2(T,a)$.

\begin{figure}
\includegraphics[scale=0.5]{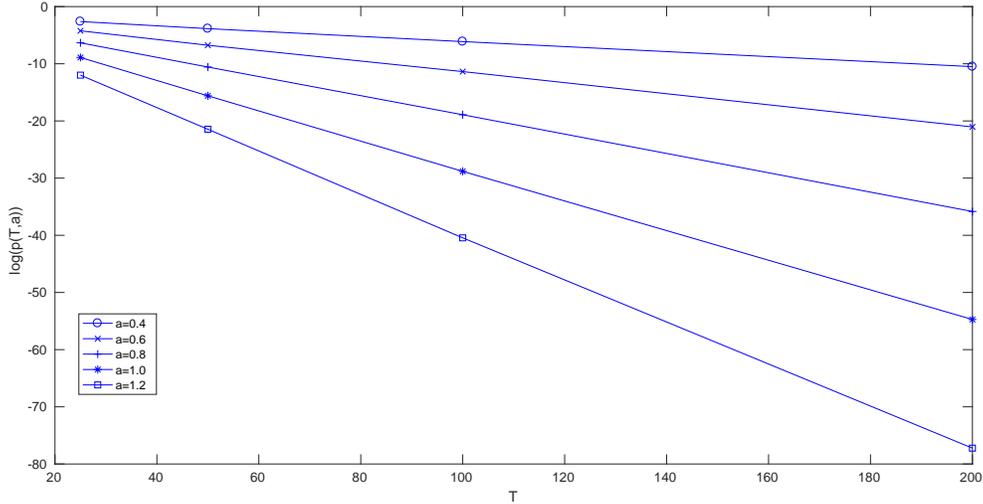}
\caption{Evolution of $\log\bigl(\hat p(T,a)\bigr)$ as a function of $T$, for different values of $a\in\left\{0.4,0.6,0.8,1.0,1.2\right\}$, see Table~\ref{tab:Ou_averages}.}
\label{fig:averageOU}
\end{figure}

In the large time limit $T\to\infty$, since the law of $Y(T)$
converges to a centered Gaussian with variance $1$, $Y(T)$ satisfies a large deviation
principle, with rate function $I$ defined by: for all $a>0$,
\[
\underset{T\to\infty}\lim~-\frac{1}{T}\log\Bigl(\mathbb{P}\bigl(Y(T)>a\bigr)\Bigr)=I(a)=\frac{a^2}{4}.
\]
In the numerical experiment, we illustrate the potential of the AMS
algorithm to estimate the large deviations rate function. Notice that
in the large $T$ limit, the probability is extremely small and in practice cannot be estimated by the vanilla splitting strategy. The
estimate of the rate function $\hat I(a)$ is obtained by a regression procedure, see
Figure~\ref{fig:averageOU}. In addition to statistical error, two sources of numerical error are identified: values of $T$ may not be sufficiently large, and the discretization of the dynamics and of the computation of temporal averages introduces a bias. The results, reported in Table~\ref{tab:Ou_averages}, show the interest of this approach to estimate large deviations rate functionals.

\begin{table}[htbp]
\begin{center}
\begin{tabular}{||c||c||c||c||c||c|c||}
\hline
$a$ & $\hat{p}(T=25,a)$ & $\hat{p}(T=50,a)$ & $\hat{p}(T=100,a)$ &
                                                                   $\hat{p}(T=200,a)$
  &  $\hat I(a)$ & $\frac{a^2}{4}$\\
 & $\hat{\sigma}^2(T=25,a)$ & $\hat{\sigma}^2(T=50,a)$ & $\hat{\sigma}^2(T=100,a)$ & $\hat{\sigma}^2(T=200,a)$ &  & \\
\hline
$0.4$ & $7.28~10^{-2}$ & $2.12~10^{-2}$ & $2.22~10^{-3}$ & $2.75~10^{-5}$ & $0.045$ & $0.040$\\
$-$ & $2.02~10^{-5}$ & $3.25~10^{-6}$ & $1.27~10^{-7}$ & $3.85~10^{-11}$ & $-$ & $-$\\
\hline
$0.6$ & $1.45~10^{-2}$ & $1.16~10^{-3}$ & $1.16~10^{-5}$ & $7.28~10^{-10}$ & $0.096$ & $0.090$\\
$-$ & $1.36~10^{-6}$ & $5.85~10^{-8}$ & $2.21~10^{-10}$ & $6.44~10^{-19}$ & $-$ &$-$\\
\hline
$0.8$ & $1.76~10^{-3}$ & $2.57~10^{-5}$ & $6.12~10^{-9}$ & $2.73~10^{-16}$ & $0.169$ & $0.160$\\
$-$ & $8.01~10^{-8}$ & $3.32~10^{-10}$ & $1.44~10^{-16}$ & $1.98~10^{-31}$ & $-$ & $-$\\
\hline
$1.0$ & $1.37~10^{-4}$ & $1.67~10^{-7}$ & $3.06~10^{-13}$ & $1.71~10^{-24}$ & $0.261$ & $0.250$\\
$-$ & $1.47~10^{-9}$ & $2.07~10^{-14}$ & $1.16~10^{-25}$ & $9.32~10^{-48}$ & $-$ & $-$\\
\hline
$1.2$ & $6.21~10^{-6}$ & $4.83~10^{-10}$ & $2.71~10^{-18}$ & $2.88~10^{-34}$ & $0.373$ & $0.360$\\
$-$ & $1.40~10^{-11}$ & $6.69~10^{-19}$ & $3.49~10^{-35}$ & $5.21~10^{-67}$ & $-$ & $-$\\
\hline
\end{tabular}
\caption{Temporal averages for an Ornstein-Uhlenbeck process. $n_{\rm rep}=10^3$ and $M=10^2$}
\label{tab:Ou_averages}
\end{center}
\end{table}

\subsection{Lorenz model}

We consider the following stochastic version of the $3$-dimensional
Lorenz system, see~\cite{BeckZuev:17,Lorenz:63} for similar numerical experiments:
\[
\begin{cases}
dX_1^\beta(t)=\sigma\bigl(X_2(t)-X_1(t)\bigr)dt+\sqrt{2\beta^{-1}}dW(t),\\
dX_2^\beta(t)=\bigl(rX_1(t)-X_2(t)-X_1(t)X_3(t)\bigr)dt,\\
dX_3^\beta(t)=X_1(t)X_2(t)-bX_3(t),
\end{cases}
\]
which depends on parameters $\sigma$, $r$, $b$ and $\beta$. The parameters are given the following values in this section: $\sigma=3$, $r=26$ and $b=1$.

Consider first the deterministic case, {\it i.e.} $\beta=\infty$. Then the system admits three unstable equilibria, and one of them is 
\[
x^\star=\bigl(\sqrt{b(r-1)},\sqrt{b(r-1)},r-1\bigr)=(5,5,25).
\]
Let the initial condition be given by
$X^\infty(0)=x^\star+\frac{1}{2}(1,1,1)$. Then, one has the following
stability result~\cite{BeckZuev:17}: for all $t\ge 0$, $\Phi\bigl(X^\infty(t)\bigr)\le
1$, where
\[
\Phi(x)=\frac{x_1^2}{(r+\sigma)^2\frac{b}{\sigma}}+\frac{x_2^2}{(r+\sigma)^2 b}+\frac{\bigl(x_3-(r+\sigma)\bigr)^2}{(r+\sigma)^2}.
\]

When noise is introduced in the system, {\it i.e.} $\beta<\infty$, we are interested in the estimation of the probability
\[
\mathbb{P}\bigl(\Phi(X^\beta(T))>1\bigr),
\]
with threshold $a=1$.

In the numerical experiments, $\sqrt{2\beta^{-1}}=3$, and the discretization is performed using the explicit Euler-Maruyama method, with time-step size $\Delta t=10^{-2}$. The sample size is $M=10^4$, and the number of replicas is $n_{\rm rep}=10^3$.

\begin{table}[htbp]
\begin{center}
\begin{tabular}{||c||c|c||c|c||}
\hline
$T$ & $\hat{p}$ & confidence interval & $\hat{\sigma}^2$ & $\frac{-\hat{p}^2\log(\hat{p})}{n_{\rm rep}}$\\
\hline
$5$ & $1.413~10^{-5}$ & $[1.388~10^{-5},1.438~10^{-5}]$ & $1.648~10^{-10}$ & $2.230~10^{-12}$ \\
$10$ & $2.607~10^{-5}$ & $[2.534~10^{-5},2.681~10^{-5}]$ & $1.409~10^{-9}$ & $7.174~10^{-12}$\\
$15$ & $2.709~10^{-5}$ & $[2.609~10^{-5},2.809~10^{-5}]$ & $2.592~10^{-9}$ & $7.718~10^{-12}$\\
$20$ & $2.594~10^{-5}$ & $[2.484~10^{-5},2.704~10^{-5}]$ & $3.158~10^{-9}$ & $7.106~10^{-12}$\\
\hline
\end{tabular}
\caption{Lorenz model. $n_{\rm rep}=10^3$ and $M=10^4$}
\label{tab:Lorenz}
\end{center}
\end{table}

This numerical experiment thus illustrates the potential of the adaptive multilevel splitting algorithms introduced in this article, for applications to complex, nonlinear, stochastic models.

\subsection{Driven periodic diffusion}

We finally consider an example taken from~\cite{FerreTouchette:18}. In this section, we consider the SDE on the unit circle, {\it i.e.} on the torus $\mathbb{T}$,
\[
dX(t)=\bigl(-V'(X(t))+\gamma\bigr)dt+\sqrt{2}dW(t),
\]
where the potential energy function $V(x)=\cos(2\pi x)$ is periodic,
and $\gamma\in \mathbb{R}$. If $\gamma\neq 0$, this is called a
non-equilibrium process since the drift term $-V'(x) +\gamma$ is not
the derivative of a function defined on the torus $\mathbb{T}$. In the remainder of this section, let $\gamma=1$, and let the initial condition in the simulation be $X_0=0$. The discretization is performed using the Euler-Maruyama method, with time-step size $\Delta t=10^{-2}$.

We are interested in the behavior of $\frac{X_T}{T}$, when $T\to \infty$, more precisely we apply the AMS algorithm to estimate
\[
p(T,a)=\mathbb{P}\bigl(X(T)>aT\bigr)=\mathbb{P}\left(\frac{X(T)}{T}>a\right).
\]
Following the same approach as for the temporal averages of the Ornstein-Uhlenbeck process, a large deviations rate function
\[
I(a)=\underset{T\to\infty}\lim~-\frac{1}{T}\log\bigl(p(T,a)\bigr),
\]
is estimated, based on estimators of the probability $p(T,a)$ for several values of $T$.

In this numerical experiment, we compare two ways of applying the AMS algorithm, with the new score function $\xi_{\rm new}$ but with different processes: considering either the process $\bigl(X(t)\bigr)_{0\le t\le T}$ with the threshold $aT$, or the process $\bigl(Y(t)=\frac{X(t)}{t}\bigr)_{0<t\le T}$, with the threshold $a$. Numerical values for different choices of $a$, $T$ and $\nrep$, of the associated estimators $\hat{p}^X(T,a)$ and $\hat{p}^Y(T,a)$, and of the empirical variances $\hat{\sigma}^{2,X}(T,a)$ and $\hat{\sigma}^{2,Y}(T,a)$ are reported in Table~\ref{tab:periodicdriven} below.

It is observed that $\hat{\sigma}^{2,Y}(T,a)<\hat{\sigma}^{2,X}(T,a)$, but a fair comparison requires to take into account the (average) computational cost. Thus the relative efficiency ${\rm Eff}(Y|X)$ of using the process $Y$ instead of $X$, is computed as the ratio
\[
{\rm Eff}(Y|X)=\frac{\hat{\sigma}^{2,X}{\rm comp. time}(X)}{\hat{\sigma}^{2,Y}{\rm comp. time}(Y)},
\]
where $\frac{{\rm comp. time}(X)}{{\rm comp. time}(Y)}$ is the ratio
of the total computational times for the experiments using $X$ and $Y$
respectively. The values of ${\rm Eff}(Y|X)$ in this numerical
experiment are reported in Table~\ref{tab:periodicdriven}. We observe
that ${\rm Eff}(Y|X)>1$ which means that the algorithm is more
efficient using the process $Y$ than the process $X$. To have a comparison with the committor score function, since the value of $p(T,a)$ is not known, an approximation of the optimal variance is computed using the estimator $\hat{p}^Y(T,a)$.

Estimators $\hat{I}(a)$ of the large deviations rate function $I(a)$ are estimated by a regression procedure (with respect to $T$) using the estimators $\hat{p}^{Y}(T,a)$, for several values of $a$. The numerical values are in excellent agreement with the numerical experiments in~\cite{FerreTouchette:18}. The AMS algorithm can thus be an efficient tool to estimate large deviations rate functions.

{\small
\begin{table}[htbp]
\begin{center}
\begin{tabular}{||c|c|c||c|c||c|c|c|c||c||}
\hline
$a$ & $T$ & $n_{\rm rep}$ & $\hat{p}^X(T,a)$ & $\hat{p}^Y(T,a)$ & $\hat{\sigma}^{2,X}(T,a)$ & $\hat{\sigma}^{2,Y}(T,a)$ & $\frac{-(\hat{p}^Y(T,a))^2\log(\hat{p}^Y(T,a))}{n_{\rm rep}}$ & ${\rm Eff}(Y|X)$ & $\hat{I}(a)$ \\
\hline
$0.8$ & $100$ & $10^2$ & $8.483~10^{-2}$ & $8.489~10^{-2}$ & $7.136~10^{-4}$ & $2.487~10^{-4}$ & $1.777~10^{-4}$ & $1.0$ & \\
$-$ & $200$ & $-$ & $2.647~10^{-2}$ & $2.776~10^{-2}$ & $1.832~10^{-4}$ & $3.519~10^{-5}$ & $2.762~10^{-5}$ & $1.7$ & $0.0112$\\
\hline
$1$ & $50$ & $10^3$ & $1.529~10^{-2}$ & $1.505~10^{-2}$ & $7.046~10^{-6}$ & $1.513~10^{-6}$ & $9.504~10^{-7}$ & $3.0$ & \\
$-$ & $100$ & $-$ & $1.026~10^{-3}$ & $1.085~10^{-3}$ & $2.586~10^{-7}$ & $3.879~10^{-8}$ & $8.036~10^{-9}$ & $2.8$ & $0.0526$\\
\hline
$1.25$ & $50$ & $10^3$ &  $1.374~10^{-4}$ & $1.311~10^{-4}$ & $1.227~10^{-8}$ & $1.355~10^{-9}$ & $1.537~10^{-10}$ & $4.9$ &\\
$-$ & $100$ & $-$ & $8.941~10^{-7}$ & $1.017~10^{-7}$ & $1.048~10^{-13}$ & $1.585~10^{-15}$ & $1.665~10^{-16}$ & $35$ & $0.189$\\
\hline
\end{tabular}
\caption{Estimates of $\mathbb{P}(X_T>aT)$ and of $I(a)$ for the
  Periodic driven diffusion. The sample size is $M=100$.}
\label{tab:periodicdriven}
\end{center}
\end{table}
}

\section*{acknowledgments}
The authors would like to thank G.~Ferr\'e and H.~Touchette for
discussions concerning the numerical experiments, and F.~C\'erou,
A.~Guyader and M.~Rousset for stimulating discussions at early stages
of this work. The work of T. Leli\`evre is supported by the European Research Council under the European Union's Seventh Framework Programme (FP/2007-2013) / ERC Grant Agreement number 614492.
\section*{acknowledgments}

\bibliographystyle{abbrv}

\end{document}